\documentclass[12pt]{article}
\usepackage{amssymb}
\usepackage[dvipdfm]{graphicx}
\usepackage{amsfonts}
\usepackage{latexsym}
\usepackage{amsthm}
\usepackage{amsmath}

\usepackage{amssymb, amscd}

\usepackage{url}

\newtheorem{problem}{Problem}
\newtheorem{theo}[problem]{Theorem}
\newtheorem{rem}[problem]{Remark}

\newtheorem{defin}[problem]{Definition}
\newtheorem{prop}[problem]{Proposition}
\newtheorem{cor}[problem]{Corollary}

\newtheorem{exam}[problem]{Example}

\begin{document}

 \title{{Hyperplane mass equipartition problem \\ and the shielding functions of Ramos\footnote{This paper
 is an updated and expanded version of \cite{Vamos-1}.}}}

\author{Sini\v sa T.\ Vre\' cica\\ {\small Faculty of Mathematics}\\[-2mm] {\small University of Belgrade}
\\[-2mm] {\small vrecica$@$matf.bg.ac.rs}
 \and Rade  T.\ \v Zivaljevi\' c\\ {\small Mathematical Institute}\\[-2mm] {\small SASA, Belgrade}\\[-2mm]
 {\small rade$@$mi.sanu.ac.rs} }
\date{March 1, 2016}

\maketitle 
\begin{abstract}
We give a proof (based on methods and ideas developed in
\cite{Guide2, MVZ, Z08}) of the result of Ramos \cite{Ram} which
claims that two finite, continuous Borel measures $\mu_1$ and
$\mu_2$ defined on $\mathbb{R}^5$ admit an equipartition by a
collection of three hyperplanes.  Our proof illuminates one of the
central methods developed and used in our earlier papers and may
serve as a good `test case' for addressing (and resolving) the
`issues' raised in \cite{BFHZ-1}; see Sections~\ref{sec:intro} and
\ref{sec:remarks} for an outline and summary. We also offer a
degree-theoretic interpretation of the `parity calculation method'
developed in \cite{Ram} and demonstrate that, up to minor
corrections or modifications, it remains a rigorous and powerful
tool for proving results about  mass equipartitions.
\end{abstract}

\section{Introduction}\label{sec:intro}

The Gr\"{u}nbaum-Hadwiger-Ramos hyperplane mass equipartition
problem \cite{Gru, Hadw, Avis, Ram, Guide2, MVZ, Z08, Ziv2015,
BFHZ-1, BFHZ-2} has been for decades one of the important test
problems for applications of topological methods in discrete
geometry.

\medskip
The problem came into the focus again with the appearance of the {
`critical review'} \cite{BFHZ-1} which, as claimed by the authors,
included the { `documentation of essential gaps'} in the proofs of
some of the earlier papers. In turn this led to an interesting and
important academic discussion about the validity, scope and
applicability of the previously used methods.

\medskip
In this paper we address the central objections raised in
\cite{BFHZ-1} (the reader will find a brief summary in
Section~\ref{sec:criticism} and concluding comments in
Section~\ref{sec:remarks}).

\medskip

We begin with a  proof of the result of Edgar Ramos \cite{Ram}
addressing the problem of equipartition of two finite, continuous
Borel measures $\mu_1$ and $\mu_2$ defined on $\mathbb{R}^5$ by a
collection of three hyperplanes. The central idea of the proof
(the `moment curve based evaluation of the topological
obstruction') originally appeared in one of our papers almost
twenty years ago, see \cite[Proposition~4.9]{Guide2}.

\medskip

As a corollary we obtain a new proof of the result
$\Delta(1,4)\leq 5$ (also due to Ramos) which states that each
continuous measure in $\mathbb{R}^5$ admits an equipartition by
$4$ hyperplanes.

\smallskip Despite the criticism and doubts raised in
\cite{BFHZ-1}\footnote{According to the `critical review'
\cite[Table 2]{BFHZ-1}, the inequality $5\leq \Delta(2,3)\leq 8$
was the only available information about the number $\Delta(2,3)$
at the time when the preprint \cite{BFHZ-1} was submitted.}, we
prefer to interpret our evaluation $\Delta(2,3)=5$ as a different
proof rather than the first complete proof of this result. As
already emphasized, our proof illuminates one of the central
methods developed in our earlier papers and we see it as a good
`laboratory test case' for discussing some of the `issues' raised
in \cite{BFHZ-1}.

\medskip
In the second half of the paper (beginning with
Section~\ref{sec:Ramos}) we give an exposition of the Ramos
`parity calculation method' \cite{Ram} emphasizing some of the key
ideas, including the concept of the {\em shielding function}
(Sections~\ref{sec:shield-principle} and
\ref{sec:shield-revisited}).

\medskip
Detailed comments and concluding remarks, summarizing our current
knowledge and opinion about the mass equipartition questions
discussed in this paper, are collected in
Section~\ref{sec:remarks}. Finally the Appendix
(Section~\ref{sec:appendix}) is a short outline of fundamental
ideas and facts about {\em transversality of equivariant maps}
which should make the paper self-contained and easier to read.

\medskip
Our ambition and the main objective in this paper was to address
all main `issues' raised by the authors of \cite{BFHZ-1}. As it
turned out there is actually only one central `issue', related to
the equivariant obstruction theory for non-free group actions. We
shall demonstrate that the `issue' disappears once we properly
interpret the role of {\em shielding functions} introduced already
by E.~Ramos precisely for this purpose, the fact well known to the
authors and many other experts in the field.

\subsection{The CS/TM-scheme and the criticism of \cite{BFHZ-1} }
\label{sec:criticism}

For the reader's orientation here we place the criticism of
\cite{BFHZ-1} in the context of the general CS/TM-scheme for the
mass equipartition problem. We hope that this outline may help the
reader understand the main objection(s) of \cite{BFHZ-1} and serve
as an introduction to the rather obvious (and well known among
specialists) remedy for the problem. In particular we emphasize
the role of the `shielding functions' which were originally
introduced by Ramos in \cite{Ram} precisely to avoid these
difficulties.

\bigskip
A configuration space (I), the associated test space and the test
map (II), and a topological result of Borsuk-Ulam type (III), are
basic ingredients of the {\em Configuration space/test map method}
(CS/TM-scheme) for applying equivariant topological methods in
discrete geometry and combinatorics, see \cite{Guide2, Z04} for an
overview.

\medskip
The general set-up for the CS/TM-scheme in the case of the  mass
equipartition problem was proposed by Ramos in \cite{Ram}. He in
particular identified proper configuration spaces and the test
maps (steps (I) and (II)), which have been without essential
change used in all subsequent publications.

\medskip
In turn this led to the general agreement that the central
difficulty in the problem is to establish the non-triviality of
the associated topological obstruction (part (III)).

\medskip
There have been two general methods to approach
(III).\footnote{The methods applied in \cite{Ziv2015} are based on
somewhat different ideas so its presentation is postponed for a
subsequent publication.}

\begin{enumerate}
 \item[(A)] The `parity count method' of Ramos \cite{Ram};
 \item[(B)] The `moment curve based evaluation of the topological
obstruction', introduced in  \cite{Guide2} and subsequently
developed in  \cite{MVZ, Z08}.
\end{enumerate}
The implementation of both of these methods was criticized in
\cite{BFHZ-1} and the authors of this paper claimed to have found
`essential gaps' in the proofs with a conclusion that {\em `the
approaches employed cannot work'} (see \cite{BFHZ-1}, the end of
the page 2).

\medskip
As the authors of some of the criticized papers, following the
dictum that one should `consistently question one's own findings',
we took these claims very seriously. Moreover,  our professional
curiosity was aroused and we wanted to understand the deeper
nature of these claims.

\medskip
Here is the summary of our response (more details can be found in
Section~\ref{sec:remarks} and elsewhere in the paper). We found
that the criticism of \cite{BFHZ-1} really applies to the `test
map phase' (step (II)) of the CS/TM-procedure and that it can be
summarized as follows. The action of the symmetry group on the
configuration space is not free. If the closed subspace of all
singular orbits is removed, one obtains a truncated space (open
manifold) where the topological obstruction is almost certainly
equal to zero (and therefore an equivariant map should exist).

\medskip
For illustration there certainly exists a
$\mathbb{Z}_2$-equivariant map $f : S^3\setminus\{a,
-a\}\rightarrow \mathbb{R}^3$ without zeros, i.e.\ the Borsuk-Ulam
theorem is no longer true if one removes two antipodal points from
the domain. However, if the map $f$ can be extended to a
$\mathbb{Z}_2$-equivariant map $g : S^3\rightarrow \mathbb{R}^3$
without zeros in $\{a,-a\}$ (or if it is equivariantly homotopic
to such a map) than the Borsuk-Ulam theorem holds for $f$ as well.

\medskip
In practise $f$ is already defined on the whole configuration
space and it is guaranteed that it has no zeros in the singular
set by the {\em shielding functions} (called the `shield
functions' by Ramos in \cite{Ram}).

\medskip
E.\ Ramos was fully aware of this technical difficulty and the
shielding functions were introduced by him in \cite{Ram} precisely
for this purpose. The reader is referred to
Section~\ref{sec:shield-principle} and
Section~\ref{sec:shield-revisited} for a more detailed explanation
of the importance of shielding functions and their role in methods
(A) and (B).

\medskip
In summary, with this clarification, both the methods (A) and (B)
are fully applicable and the proofs and results obtained by their
application are correct.


\section{Mass equipartitions by hyperplanes}
\label{sec:mass}

The reader is referred to \cite{Ram, Guide2, MVZ, Ziv2015, BFHZ-1}
for an overview of known results and the history of the general
measure equipartition problem by hyperplanes. Recall that the
problem has its origins in the papers of Gr\"{u}nbaum \cite{Gru}
and Hadwiger \cite{Hadw}, with the papers of Steinhaus \cite{Ste}
and Stone and Tukey \cite{Sto-Tuk} as important predecessors.

\subsection{CS/TM-scheme for the mass equipartition problem}
\label{sec:set-up}

A collection $\mathcal{A}=\{A_1,\ldots, A_j\}$ of Lebesgue
measurable sets in $\mathbb{R}^d$ admits an {\em equipartition} by
a collection $\mathcal{H}=\{h_1,\ldots, h_k\}$ of hyperplanes if
$m(A_i\cap O) = (1/2^k)m(A_i)$ for each $i=1,\ldots, j$ and each
of the $2^k$ hyperorthants $O$ associated to $\mathcal{H}$.

\medskip
More generally a collection $\mathcal{M}=\{\mu_1,\ldots, \mu_j\}$
of continuous, finite, Borel measures defined on $\mathbb{R}^d$
admits an equipartition by $\mathcal{H}$ if $\mu_i(O) =
1/2^k\mu_i(\mathbb{R}^d)$ for each $i=1,\ldots, j$.

\medskip
The `equipartition number' $\Delta(j,k)$ is defined as the minimum
dimension $d$ of the ambient space $\mathbb{R}^d$ such that each
collection $\mathcal{M}$ of $j$ continuous measures admits an
equipartition by some collection $\mathcal{H}$ of $k$ hyperplanes.
We also say that a triple $(d,j,k)$ is {\em admissible} if
$\Delta(j,k)\leq d$.

\medskip
Following \cite{Ram} and \cite{MVZ} the compactified {\em
configuration space} for the general mass equipartition problem is
the manifold $M_{d,k} = (S^d)^k$ where $h = (h_1,\ldots, h_k)\in
M_{d,k}$ is an ordered $k$-tuple of oriented hyperplanes
(including the hyperplanes `at infinity').

\medskip
Given a $0$-$1$-sequence (alternatively a $(+-)$--sequence) $J =
(j_1j_2\ldots j_k)\in 2^{[k]}$ and a $k$-tuple $h\in M_{d,k}$, the
associated half-spaces are $h_1^{j_1},\ldots, h_k^{j_k}$ and the
{\em test function} $a_J^\mu(h) := \mu(\bigcap_{\nu=1}^k
h_\nu^{j_\nu})$ measures the amount of mass $\mu$ in the
corresponding hyperorthant. Let $\hat{J} = \{\nu\in [k]\mid j_\nu
= 1\}$ be the subset of $[k]$ determined by  $J$.

\medskip
Following \cite{Ram} the collection $\{a_J^\mu\mid J\in 2^{[k]}\}$
of test functions is (via a Discrete Fourier Transform) replaced
by the functions,
\begin{equation}\label{eqn:DFT}
f_I^\mu(h) = f^\mu_{i_1\ldots i_k}(h_1,\ldots, h_k) = \sum_{J\in
2^{[k]}}(-1)^{\langle I, J\rangle}a_J^\mu(h)
\end{equation}
where $\langle I, J\rangle =  i_1j_1+\ldots + i_kj_k = \vert
\hat{I}\cap \hat{J}\vert$ is the cardinality of the set
$\hat{I}\cap \hat{J}$.

\begin{rem}{\rm
The reader should observe that in the CS/TM-scheme described here
we tacitly use the continuity properties of measures when we
extend functions $a_J^\mu$ and $f_I^\mu$ to the whole
configuration space $M_{d,k}$  (which includes hyperplanes `at
infinity'). }
\end{rem}

\subsection{ $\Delta(2,3) = 5$}
\label{sec:delta-2-3}

The equipartition problem attracted new audience and received
wider recognition with the appearance of the paper of E.~Ramos
\cite{Ram} who introduced new technique and obtained many new
results about the function $d = \Delta(j,k)$ including the
following,
\begin{equation}\label{eqn:Ramos}
4\leq \Delta(1,4)\leq 5 \qquad \Delta(3,2) = 5 \qquad 7\leq
\Delta(3,3)\leq 9 \qquad \Delta(4,2) = 6.
\end{equation}
Among the most interesting is his claim that $\Delta(2,3) = 5$
which allowed him to prove, by a simple reduction, that
$\Delta(1,4)\leq 5$.

\medskip
Here we give a different proof of a slightly more general result.
The generalization may be of some independent interest, however it
primarily exemplifies the phenomenon that a strengthened statement
may be sometimes easier to prove, cf.\
\cite[Propostition~4.9]{Guide2} for an early example in the
context of equipartitions of masses by hyperplanes.

\medskip
Recall (Section~\ref{sec:set-up}) that if  $h$ is an oriented
hyperplane then for $\epsilon\in\{+, -\}$ the associated closed
half-space is denoted by $h^\epsilon$.

\begin{theo}\label{thm:glavna}
Suppose that $\mu_1, \mu_2, \mu_3$ are three continuous, finite,
non-negative Borel measures defined on $\mathbb{R}^5$. Then there
exist three hyperplanes $h_1, h_2, h_3$ in $\mathbb{R}^5$ forming
an equipartition for measures $\mu_1$ and $\mu_2$ such that one of
them {\rm (}$h_i$ for some $i\in [3]${\rm )} is a bisector of
$\mu_3$ in the sense that $\mu_3(h_i^+) = \mu_3(h_i^-)$.
\end{theo}

\medskip\noindent
{\bf Proof:} The theorem says that in addition to being an
equipartition for the measures $\mu_1, \mu_2$ it can be always
achieved that one of the hyperplanes $h_1, h_2, h_3$ is a halving
hyperplane for a third measure $\mu_3$ (which is also prescribed
in advance).

\medskip
Following the usual {\em configuration space/test map scheme}, see
\cite{MVZ} and Section~\ref{sec:set-up}, the configuration space
(parameterizing all triples $(h_1, h_2, h_3)$ of oriented
hyperplanes in $\mathbb{R}^5$ including the hyperplanes `at
infinity') is $M = S^5\times S^5\times S^5$. The group of
symmetries acting on $M$ by permuting the hyperplanes and changing
their orientation is the group $G = \mathbb{Z}_2^{\oplus 3}\rtimes
S_3$. This group arises also as the group of symmetries of the
$3$-dimensional cube.

\medskip
The real regular representation $\mathbb{R}[\mathbb{Z}_2^{\oplus
3}]$ of $\mathbb{Z}_2^{\oplus 3}$ is also a $G$-representation.
The test space $V_i$ for each of the measures $\mu_i, \, i=1,2$ is
the $G$-representation of dimension $7$ arising by subtracting
from $\mathbb{R}[\mathbb{Z}_2^{\oplus 3}]$ the trivial
$1$-dimensional representation,
$$
V_i = \{\sum_{\epsilon} \alpha_\epsilon\cdot\epsilon\in
\mathbb{R}[\mathbb{Z}_2^{\oplus 3}] : \sum_{\epsilon}
\alpha_\epsilon = 0\}.
$$
We also need a copy of $\mathbb{R}$ to serve as the
target space for testing if one of the hyperplanes is a bisector
for $\mu_3$, so the total test space is the $G$-representation $V
= V_1\oplus V_2\oplus \mathbb{R}$.

\medskip In agreement with (\ref{eqn:DFT}) the associated `test
map',
 \begin{equation}\label{eqn:test-map}
f = (f_1, f_2, f_3) : S^5\times S^5\times S^5 \rightarrow
V_1\oplus V_2\oplus \mathbb{R}
 \end{equation}
is described as follows. Since the hyperplanes $h_i$ are oriented
each hyperorthant $O_{\epsilon} = O_{(\epsilon_1, \epsilon_2,
\epsilon_3)} = h_1^{\epsilon_1}\cap h_2^{\epsilon_2}\cap
h_3^{\epsilon_3}$ is associated an element $\epsilon =
(\epsilon_1, \epsilon_2, \epsilon_3)$ of the group
$\mathbb{Z}_2^{\oplus 3}$. By definition (for $i=1,2$),
\begin{equation}\label{eqn:test-map-bis}
f_i(h_1, h_2, h_3) = \sum_{\epsilon} [\mu_i(O_\epsilon)-
\frac{1}{8}\mu_i(\mathbb{R}^5)]\cdot\epsilon\in V_i\subset
\mathbb{R}[\mathbb{Z}_2^{\oplus 3}]
\end{equation}
The map $f_3 : S^5\times S^5\times S^5 \rightarrow \mathbb{R}$ is
defined by,
\begin{equation}\label{eqn:test_map_3}
f_3(h_1, h_2, h_3) =
(\mu_3(h_1^+)-\mu_3(h_1^-))(\mu_3(h_2^+)-\mu_3(h_2^-))(\mu_3(h_3^+)-\mu_3(h_3^-))
\end{equation}
By construction a triple $h = (h_1, h_2, h_3)$ satisfies the
conditions of the theorem if and only if $h$ is a zero of the test
map (\ref{eqn:test-map}). The test map is clearly $G$-equivariant
so it is sufficient to show that there does not exist a
$G$-equivariant map $f : M\rightarrow V\setminus\{0\}$ (arising
from measures, as in the construction above).

\begin{rem}\label{rem:uslov}
{\rm  If $\mu_1(\mathbb{R}^5) = \mu_2(\mathbb{R}^5) = 0$ then
there is a `trivial' solution  of the equation $f(h)=0$. Indeed,
if $H$ is a bisector of both $\mu_1$ and $\mu_2$ then $h = (H, H,
H)$ is a clearly an equipartition for $\{ \mu_1, \mu_2\}$. For
this reason we tacitly assume (throughout the paper) that either
$\mu_1(\mathbb{R}^d)>0$ or $\mu_2(\mathbb{R}^d)>0$ (or both). }
\end{rem}

\subsection{Shielding functions}
\label{sec:shield-principle}

The action of $G$ on $M = S^5\times S^5\times S^5$ is not free.
The singular set \[S = \bigcup_{x,y\in S^5}( G\cdot (x,x,y)\cup
G\cdot (x,-x,y))\] consists of ordered triples having two equal or
two antipodal elements. In light of Remark~\ref{rem:uslov} we are
allowed to assume  that $\mu_i(\mathbb{R}^5) = \lambda > 0$ for at
least one $i\in \{1,2\}$, say for $i=1$. For most of the
interesting measures (including positive measures) both of these
two conditions are satisfied.

\medskip
As an immediate consequence we deduce that $f(S)\subset
V\setminus\{0\}$. Moreover, all maps $f = f_{(\mu_1,\mu_2,\mu_3)}
: S \rightarrow V\setminus\{0\}$ arising from measures such that
$\mu_1(\mathbb{R}^5)>0$ {\bf are linearly $G$-homotopic}. Indeed,
the value of $f_1(x,x,y)\in \mathbb{R}[\mathbb{Z}_2^{\oplus 3}]$
has a non-zero coefficient $-(1/8)\lambda$ associated to the
element $\epsilon = (+,-,\epsilon_3)$ while the value of
$f_i(x,-x,y)$ has the same non-zero coefficient at $\epsilon =
(+,+,\epsilon_3)$. This is (in light of the definition of the test
map (\ref{eqn:test-map-bis})) a simple consequence of the fact
that whenever a hyperorthant is empty the corresponding
coefficient is equal to $-(1/8)\mu_i(\mathbb{R}^5)\neq 0$.

\medskip
Let $g: S \rightarrow V\setminus\{0\}$ be a representative of this
$G$-homotopy class, for example $g$ can be chosen to be the
restriction of a map $f_{\nu_1,\nu_2,\nu_3}$ for a particular
choice of measures $\nu_1, \nu_2, \nu_3$. Theorem~\ref{thm:glavna}
will be proved if we show that the first (and only) obstruction
class $\omega = \omega_g$ for extending equivariantly the map $g$
to $M = S^5\times S^5\times S^5$ is non-zero.

\begin{defin}\label{def:shield}{\rm(Shielding functions)}
 The functions that keep the zeros of the test map away from
the singular set (where the action of the group is not free) are
following Ramos {\rm \cite{Ram}} called the {\em shielding or
shield functions}. The fact that all test maps arising from
measures are non-zero and $G$-homotopic on the singular set is
referred to as the {\em `shielding functions homotopy principle'}.
\end{defin}

\subsection{Calculation of the obstruction $\omega$}
\label{sec:obstr-omega}

The obstruction $\omega$ described in the previous section lives
in the relative equivariant cohomology group $H_G^{15}(M, S; W)$
where $W$ is the $G$-module $\pi_{14}(V\setminus\{0\})$. By taking
a small $G$-equivariant open tubular neighborhood $U$ of $S$ we
observe that $\omega$ can be evaluated in the group $H_G^{15}(N,
\partial N; W)$ where $N = M\setminus U$ is a compact $G$-manifold
with boundary.

\medskip
In light of the equivariant Poincar\'{e}-Lefschetz duality this
group is isomorphic to the equivariant homology group $H^G_0(N,
W\otimes \pi)$ where $\pi$ is the orientation character describing
the action of $G$ on the (relative) fundamental class of $N$. This
group is isomorphic to one of the groups $\mathbb{Z}_2$ or
$\mathbb{Z}$ and the corresponding dual of $\omega$ can be
evaluated in this group by a geometric argument.

\medskip
Note that since the action of $G$ on the manifold with boundary
$N$ is free, by passing to the quotient manifold $N/G$ we can
actually use the usual version of Poincar\'{e}-Lefschetz duality
(with local coefficients). However there is a shortcut which
allows us to complement (bypass) the homological algebra related
to the (equivariant) Poincar\'{e}-Lefschetz duality and replace it
by a direct geometric argument.\footnote{Here we follow the same
general strategy applied in \cite{MVZ}, both in the overall
technical set-up and in the computational details (which are of
course much more complex in the case $(8,5,2)$ studied in the
paper \cite{MVZ}).}

\begin{prop}\label{prop:claim}
Suppose that $f : M \rightarrow V$ is a $G$-equivariant extension
of $g : S \rightarrow V\setminus\{0\}$. Moreover, assume that $f$
is smooth outside of a small tubular neighborhood $U$ of $S$ {\rm
(}$g(U)\subset V\setminus\{0\}${\rm )} and that $f$ is transverse
to $0\in V$. The set $f^{-1}(0)$ is finite and $G$-invariant. The
number of $G$-orbits $n_f = \vert f^{-1}(0)/ G\vert$ depends on
$f$, however the parity of this number $\theta = \theta_f =_{mod\,
2} n_f\in \mathbb{Z}_2$ is the same for all extensions $f$ of $g$.
In particular if this number is odd then each $G$-extension of $g$
must have a zero.
\end{prop}

\medskip\noindent
{\bf Proof:} The result is an immediate consequence of
Proposition~\ref{prop:central}.  \hfill $\square$

\begin{figure}[hbt]
\centering
\includegraphics[scale=0.75]{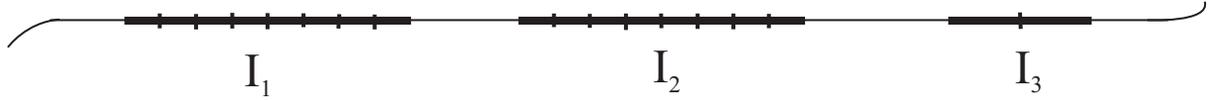}
\caption{Three hyperplanes intersect the moment curve in 15
points. } \label{fig:vamos-1}
\end{figure}

\medskip
In order to compute the obstruction $\theta\in \mathbb{Z}_2$ we
use the map $f : M \rightarrow V$ which arises from the following
choice of measures (measurable sets). Let $\Gamma$ be the moment
curve in $\mathbb{R}^5$  defined as the image $\phi(\mathbb{R})$
of the map $\phi : \mathbb{R}\rightarrow \mathbb{R}^5, \, t\mapsto
(t,t^2, \ldots t^5)$. Suppose that $I_1, I_2, I_3$ are three
disjoint, consecutive intervals on this curve
(Figure~\ref{fig:vamos-1}) and let $\nu_i$ be the measures on
$\mathbb{R}^5$ defined by $\nu_i(A) = m(\phi^{-1}(A\cap I_i))$.

\medskip
By construction $h = (h_1, h_2, h_3)\in f^{-1}(0)$ is the set of
triples of oriented hyperplanes in $\mathbb{R}^5$ such that $h$ is
an equipartition for both $\nu_1$ and $\nu_2$ and in addition one
of the hyperplanes $h_i$ is a bisector of $\nu_3$.  Altogether
there are at most 15 points in the set $\Gamma\cap(h_1\cup h_2\cup
h_3)$ and all fifteen are needed (Figure~\ref{fig:vamos-1}) if
$h\in f^{-1}(0)$.

\begin{figure}[hbt]
\centering
\includegraphics[scale=0.95]{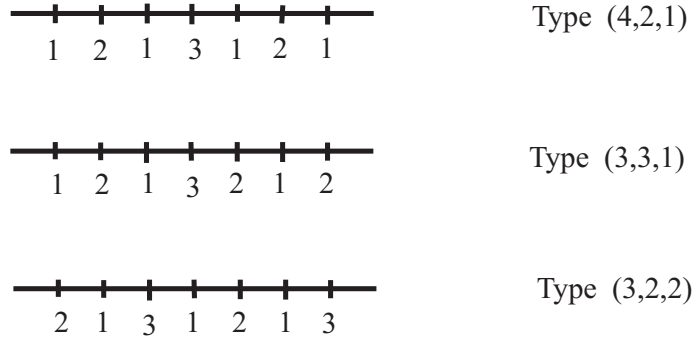}
\caption{Up to a permutation of $\{1, 2, 3\}$ there are three
possible types of equipartitions of an interval on the moment
curve by three hyperplanes.} \label{fig:vamos-2}
\end{figure}

\noindent For bookkeeping purposes let us analyze possible `types'
of elements $h = (h_1, h_2, h_3)\in f^{-1}(0)$. We begin with the
analysis which cardinalities of sets $h_i\cap I_j$ are permitted.
A closer inspection reveals that there are only three
possibilities (Figure~\ref{fig:vamos-2}) and we notice that only
some pairs of `complementary types' (associated to intervals $I_1$
and $I_2$) can appear together.

\medskip
For example one of the possibilities, abbreviated as $h = \langle
(4_1, 2_2, 1_3),(1_1, 3_2, 3_3)\rangle$ or simply as $h = \langle
(4, 2, 1),(1, 3, 3)\rangle$, describes the case where the
cardinalities of the intersections are,
\[
\vert h_1\cap I_1\vert = 4 \qquad \vert h_2\cap I_1\vert = 2
\qquad \vert h_3\cap I_1\vert = 1
\]
and
\[
\vert h_1\cap I_2\vert = 1 \qquad \vert h_2\cap I_2\vert = 3\qquad
\vert h_3\cap I_2\vert = 3.
\]

All permutations of indices are possible, for example in the same
$G$-orbit with $h$ is the element $h' = \langle(4_3, 2_1,
1_2),(1_3, 3_1, 3_2)\rangle = \langle(2,1,4), (3,3,1)\rangle$.

\medskip
The following Claim summarizes the information needed for the
evaluation of the number of $G$-orbits in $f^{-1}(0)$. Our initial
observation is that in each orbit $G\cdot h\subset f^{-1}(0)$
there is an element $h' = (h_1, h_2, h_3)$ such that the type of
the intersection $(h_1\cup h_2\cup h_3)\cap I_1$ is precisely one
of the types listed in Figure~\ref{fig:vamos-2}.

\medskip\noindent
{\bf Claim:}
\begin{enumerate}
 \item[(1)] The $I_1$-type $(4,2,1)$ can be matched with the
$I_2$-type $(1,2,4)$ in only one way, contributing {\bf 1} orbit;
 \item[(2)] The $I_1$-type $(4,2,1)$ can be matched with the
$I_2$-type $(1,3,3)$ in two ways, contributing {\bf 2} orbits;
 \item[$(3)$] The $I_1$-type $(3,3,1)$ can be uniquely matched with both
$(1,2,4)$ and $(2,1,4)$ (as the $I_2$-types), which together
contribute {\bf 2} orbits;
 \item[(4)] The $I_1$-type $(3,3,1)$
can be matched with the $I_2$-type $(2,2,3)$ in two ways,
contributing {\bf 2} orbits;
 \item[$(5)$] The $I_1$-type $(3,2,2)$ can be matched with the
$I_2$-type $(1,3,3)$ in two ways, contributing {\bf 2} orbits;
\item[(6)] The $I_1$-type $(3,2,2)$ can be matched with the
$I_2$-type $(2,3,2)$ in two ways, contributing {\bf 2} orbits;
 \item[(7)] The $I_1$-type $(3,2,2)$ can be matched with the
$I_2$-type $(2,2,3)$ in two ways, contributing {\bf 2} orbits.
 \end{enumerate}
For illustration let us check the case $(3)$. We can uniquely
choose $h = (h_1, h_2, h_3)$ in this orbit so that the $I_1$-type
is precisely the middle type shown in Figure~\ref{fig:vamos-2}.
Using this information we reconstruct the intersection of
hyperplanes $h_i$ with the interval $I_3$. If the type of this
intersection is $(1,2,4)$ the midpoint of $I_2$ belongs to $h_1$
while the intersection $h_2\cap I_2$ has two elements which also
uniquely determines their positions in $I_2$. The case when the
$I_2$-type is $(2,1,4)$ is treated similarly. The other cases
listed in the Claim are checked by a similar reasoning.

\medskip\noindent Altogether there are $13$ $G$-orbits in the set
$f^{-1}(0)$ which shows that the parity of the obstruction is
$\theta = 1$.

 \medskip
For the completion of the proof of Theorem~\ref{thm:glavna} we
should convince ourselves that all zeros $h\in f^{-1}(0)$ are
non-degenerate. This is established along the lines of the proof
of Theorem~33 in \cite{MVZ}, see also the comments on the proof of
Theorem~4 on page 291 (ibid.). The proof in our case is actually
much simpler since we are interested in the parity calculation,
i.e.\ we don't have to worry about the sign of the Jacobian
matrix.

\medskip
For a chosen $h = (h_1, h_2, h_3)\in f^{-1}(0)$ and an arbitrary
triple $l = (L_1, L_2, L_3)$ in a small neighborhood $U$ of $h$ we
observe that $l$ is determined by the fifteen numbers,
\[
\{x_1<x_2<\ldots < x_{15}\} = (L_1\cup L_2\cup L_3)\cap (I_1\cup
I_2\cup I_3)
\]
which therefore can be used as coordinating functions on $U$.

\medskip
If $I_1 = [a_1, b_1], I_2 = [a_2, b_2]$ and $I_3 = [a_3, b_3]$
then the hyperplanes $L_1, L_2, L_3$ divide the interval $I_1$ in
hyperorthants $[a_1, x_1], [x_1, x_2], \ldots, [x_7,b_1]$ and the
interval $I_2$ in hyperorthants $[a_2, x_8], [x_8, x_9], \ldots,
[x_{14},b_2]$. From here we observe that the functions,
\[
x_1-a_1, x_2-x_1, \ldots , x_7-x_6  \qquad \mbox{\rm and }\qquad
x_8-a_2, x_9-x_8, \ldots , x_{14}-x_{13}
\]
can be used as coordinates on the (truncated) test space
$V_1\oplus V_2$. By an affine change of coordinates we see that
$x_1, x_2, \ldots, x_{14}$ can be used as the coordinates on the
space $V_1\oplus V_2$ as well. From here and the fact that
$\partial f_3/\partial x_{15}\neq 0$ we easily conclude that the
corresponding Jacobian matrix is non-singular. \hfill $\square$

\begin{cor}
Each continuous measure $\mu$ in $\mathbb{R}^5$ admits an
equipartition by $4$ hyperplanes. Moreover one of these
hyperplanes can be chosen to be a common bisector of $4$
measurable sets (measures) prescribed in advance, and one of the
remaining hyperplanes is also a bisector of a chosen measurable
set.
\end{cor}

\medskip\noindent
{\bf Proof:} The case $(5,1,4)$ (one measure and the equipartition
by $4$ hyperplanes) is reduced to the case $(5,2,3)$ (two measures
and the equipartition by $3$ hyperplanes) by taking a bisector $H$
of $\mu$ and applying Theorem~\ref{thm:glavna} to the two new
measures $\mu_+, \mu_-$ defined by $\mu_+(A) = \mu(H^+\cap A)$ and
$\mu_-(A) = \mu(H^-\cap A)$. By the `ham sandwich theorem'
(applied in $\mathbb{R}^5$) the hyperplane $H$ can be chosen as a
halving hyperplane of four additional measurable sets. Similarly
by Theorem~\ref{thm:glavna} one of the remaining hyperplanes can
be chosen as a bisector for a chosen measurable set. \hfill
$\square$

\subsection{Other cases of the equipartition problem}

The method applied in our proof of Theorem~\ref{thm:glavna} is not
new, indeed it has been developed in our papers and over the years
successfully applied to many cases of the equipartition problem.

\medskip
Ramos in \cite{Ram} isolated the triples $(4,1,4), (8,5,2)$ and
$(7,3,3)$ as the first cases not covered by his method. For this
reason we focused and tested our method initially in these cases.

\medskip
The admissibility of the triple $(8,5,2)$ was established in
\cite{MVZ}. At the same time we calculated (following the same
plan as in the proof of Theorem~\ref{thm:glavna}) the obstruction
in the case of the triple $(7,3,3)$. The result turned out to be
zero and this is the reason why this observation was not
published, although it was reported in our lectures and
presentations of the (much more complicated) case $(8,5,2)$.

\medskip
Essentially the same strategy was applied by the authors of
\cite{BFHZ-2} who established (as the only new result of the
paper) the admissibility of the triple $(10,4,3)$. Moreover it is
not difficult to observe that, as long as the calculations are
dependent on the `moment curve based evaluation of the topological
obstruction' (method (B)), as described in
Sections~\ref{sec:criticism} and  \ref{sec:obstr-omega}), the
`join scheme' and the `product scheme' described in
\cite[Section~1.2.]{BFHZ-2} are equivalent! In other words the
`join scheme' cannot bring anything new that is not already
provable by the `product scheme'.

\medskip
Finally, the paper \cite{Ziv2015} used different methods to
establish the admissibility of the triple $(6\cdot 2^\nu +2,
4\cdot 2^\nu+1, 2)$ for each $\nu\geq 0$, which includes $(8,5,2)$
as a special case. This paper also came under criticism of
\cite{BFHZ-1} and these objections, being of somewhat different
nature, will be addressed in our subsequent publication.

\section{The method of Ramos as developed in \cite{Ram}}
\label{sec:Ramos}

The paper of Edgar Ramos \cite{Ram} introduced important new ideas
into the mass equipartition problem, leading to the algorithms for
calculating relevant topological obstructions. Indeed, this paper
has been for years an inspiration for all subsequent work in this
area.

\medskip
The authors of \cite{BFHZ-1} acknowledge the importance of
\cite{Ram}, however they claim that the proofs of central results
of this paper have essential gaps. If proved correct, the critical
analysis from \cite{BFHZ-1} would render obsolete not only the
proofs of the results from \cite{Ram} but the ideas and the
methods would be also affected. This in particular applies to the
`parity count' formulas from \cite{Ram} which would be probably
avoided in the future, if the criticism from \cite{BFHZ-1} is
taken for granted.

\medskip
This would be a shame since these methods (and proofs) are
(essentially) correct in a strong sense of the word. They continue
to be a valuable tool for tackling problems in this and related
areas of {\em Applied and computational algebraic topology}.

\medskip
In the following section we revisit, reprove and give a slightly
different interpretation to the parity count lemmas and formulas
from \cite{Ram}. In particular we explain why the `counterexample'
\cite[Section~7]{BFHZ-1} is not properly addressing Lemma~6.2 from
\cite{Ram} which remains correct and applicable.

\subsection{The parity count formulas from \cite{Ram}}
E. Ramos used in \cite{Ram} intricate parity count calculation (in
the setting of piecewise linear topology) to evaluate the
obstruction to the existence of equipartitions of masses by
hyperplanes. His approach was  critically analyzed in
\cite{BFHZ-1} and some of his results were doubted, in particular
the authors of \cite{BFHZ-1} questioned the validity of his proof
of the equality $\Delta(2,3) = 5$.

\medskip
We begin this exposition by observing that the parity count
formulas of Ramos \cite[p.~150]{Ram}, designed for evaluating the
`parity invariants' $P(r; X)$ and $P^+(r', \underline{r}''; X)$,
can be naturally interpreted as formulas/algorithms for evaluating
the degrees of associated maps.

\begin{defin}\label{def:Ram-degree-1}
Suppose that $r : (X, \partial(X)) \rightarrow (\mathbb{R}^n,
\mathbb{R}^n\setminus\{0\})$ is a map where $X$ is a compact,
$n$-dimensional manifold with boundary $\partial(X)$. Then by
definition,
\begin{equation}\label{eqn:deg-1}
     P(r; X) :=  {\rm deg}\{H_n(X, \partial(X); \mathbb{Z}_2)\stackrel{r_\ast}{\longrightarrow}
     H_n(\mathbb{R}^n, \mathbb{R}^n\setminus\{0\}; \mathbb{Z}_2)
     \}
\end{equation}
is the ${\rm mod}_2$-degree of the map $r$.
\end{defin}
If $X$ is a smooth or triangulated manifold such that $0$ is not a
critical value of $r$ (i.e.\ $r$ is transverse to zero $r
\pitchfork \{0\}$), then $P(r; X)$ is indeed the parity (the ${\rm
mod}_2$-cardinality) of the (finite) set $r^{-1}(0)$.

\begin{defin}\label{def:Ram-degree-2}
Suppose that $s : Y \longrightarrow \mathbb{R}^n\setminus\{0\}$ is
a map defined on a compact, $(n-1)$-dimensional manifold without
boundary. Then by definition,
\begin{equation}\label{eqn:deg-2}
     P(s; Y) :=  {\rm deg}\{H_{n-1}(Y; \mathbb{Z}_2)\stackrel{s_\ast}{\rightarrow}
     H_n(\mathbb{R}^n\setminus\{0\}; \mathbb{Z}_2)\}.
\end{equation}
\end{defin}

\begin{prop}\label{prop:drina}
Choose $v\in \mathbb{R}^n\setminus\{0\}$ and let $L = \{\lambda v
\mid \lambda\geq 0\}$ be the associated half-ray in
$\mathbb{R}^n$. Suppose that $s : Y \longrightarrow
\mathbb{R}^n\setminus\{0\}$ is a map defined on a smooth
(alternatively PL-triangulated) compact, $(n-1)$-dimensional
manifold without boundary. Assume that $s$ is transverse to the
ray $L$, $s \pitchfork L$, and let $P(s,L; Y)$ be the ${\rm
mod}_2$-cardinality of the set $s^{-1}(L)$. Then,
\begin{equation}\label{eqn:ray}
P(s; Y)  =  P(s,L; Y).
\end{equation}
\end{prop}

\medskip\noindent
{\bf Proof:} Let $\psi : \mathbb{R}^n\setminus\{0\} \rightarrow
S^{n-1}$ be the radial projection, $v \mapsto v/\| v \|$. Then $s
\pitchfork L$ if and only if $v/\| v \|$ is a regular value of the
map $\psi\circ s$.  The result follows from the well-known fact
the ${\rm mod}_2$-degree of a map $g : Y \rightarrow S^{n-1}$ can
be calculated as the ${\rm mod}_2$-cardinality of the set
$g^{-1}(a)$ for any regular value $a\in S^{-1}$. \hfill $\square$

\begin{figure}[hbt]
\centering
\includegraphics[scale=.45]{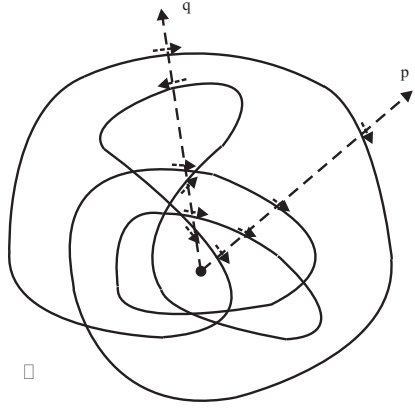}
\caption{The winding number as the number of signed intersections
with a half-ray.} \label{fig:wind}
\end{figure}
\begin{rem}{\rm The formula (\ref{eqn:ray}) in the planar case reduces to the
(${\rm mod}_2$-version) of the well-known description of the
winding number of a curve as the number of signed intersections
with a (generic) half-ray (Figure~\ref{fig:wind}). }
\end{rem}

The following `standard genericity assumption' summarizes the
conditions needed for comparison of different parity calculations
(as in \cite[Lemma~2.1.]{Ram}).

\begin{defin}\label{def:ray}{\rm (Standard genericity assumption)}
We say that a map $r = (r', r'') : (X, \partial(X))
\longrightarrow (\mathbb{R}^{n}, \mathbb{R}^n\setminus\{0\})$,
defined on a compact, $n$-dimensional manifold $X$ with boundary
$\partial(X)$, satisfies the {\bf standard genericity assumption}
if:
\begin{enumerate}
 \item[{\rm (1)}] $X$ is a smooth (alternatively PL-triangulated) manifold
such that $0$ is not a critical value of $r$.
 \item[{\rm (2)}] The restriction $s = r\vert_{\partial(X)}$ of $r$ on
 the boundary $\partial(X)$ of $X$ is transverse to $L$ where
$L=\{0\}\times \mathbb{R}^+\subset \mathbb{R}^{n-1}\times
\mathbb{R}$.
\end{enumerate}
Here $r'$ denotes the first $(n-1)$-components of $r$ and $r''$ is
the last component of $r$ while $L$ is by definition the positive
semi-axis corresponding to the last coordinate in $\mathbb{R}^n$.
\end{defin}

\begin{prop}\label{prop:standard}{\rm (\cite[Lemma~2.1.]{Ram})}
Let $r = (r', r'') : (X, \partial(X)) \rightarrow (\mathbb{R}^{n},
\mathbb{R}^n\setminus\{0\})$ be a map defined on a compact,
$n$-dimensional manifold $X$ with boundary $\partial(X)$ which
satisfies the {\em `standard genericity assumption'} in the sense
of Definition~\ref{def:ray}. Then,
\begin{equation}\label{eqn:Lema-2.1}
P(r, X) = P^+(r', \underline{r}''; \partial(X))
\end{equation}
where by definition $P^+(r', \underline{r}''; \partial(X)) = P(s,
L;
\partial(X))$.
\end{prop}

\medskip\noindent
{\bf Proof:} The equality of mapping degrees $P(r; X)$ (formula
(\ref{eqn:deg-1})) and $P(s; \partial(X))$ (formula
(\ref{eqn:deg-2})) is well-known, so the equality
(\ref{eqn:Lema-2.1}) is an immediate consequence of
Proposition~\ref{prop:drina}. \hfill $\square$

\subsection{Shielding functions revisited}
\label{sec:shield-revisited}

{\em Shielding functions} are (under the name `shield functions')
explicitly mentioned in \cite{Ram} at least three times. However
the references to these functions are ubiquitous in the paper. As
already emphasized in Section~\ref{sec:shield-principle} the role
of these functions is to {\em shield} the zeros of the
equipartition test function from appearing in the singular set $S$
where the action of the group is not free. More explicitly, in
agreement with Definition~\ref{def:shield}, $f$ is a shielding
function for $A\subset S$ if $f(x)\neq 0$ for each $x\in A$. The
following proposition, cf.\ \cite[Property 4.5.(i)]{Ram},
describes an important class of shielding functions. Recall that
for a given $0$-$1$-sequence $I\in 2^{[k]}$ the associated set is
$\hat{I} = \{\nu\in [k]\mid i_\nu = 1\}$.

\begin{prop}\label{prop:property-4.5}{\em (\cite[Property
4.5.(i)]{Ram})} Let $I = (i_1i_2\ldots i_k)\in 2^{[k]}$ and
suppose that the cardinality of the set $\hat{I}= \{\nu\in [k]\mid
i_\nu = 1\}$ is even. Suppose that $h_1,\ldots, h_k$ are oriented
hyperplanes such that $h_{i_{\nu_1}} = h_{i_{\nu_2}}$ if $\{\nu_1,
\nu_2\}\subset \hat{I}$.  Then $f_I^\mu$ is a {\em shielding
function} in the sense that,
 \begin{equation}\label{eqn:shield-2}
 f_I^\mu(h) = f^\mu_{i_1\ldots i_k}(h_1,\ldots, h_k)=1
 \end{equation}
\end{prop}

\medskip\noindent
{\bf Proof:} From the assumption that $h_{\nu_1} = h_{\nu_2}$ for
each pair $\{\nu_1, \nu_2\}\subset \hat{I}$ we deduce that
$a_J(h)$ can be non-zero only in two cases, if either $\hat{J}\cap
\hat{I}=\emptyset$ or $\hat{J}\cap \hat{I}=\hat{I}$. By definition
$f_I(h) = \sum_J~(-1)^{\vert \hat{I}\cap \hat{J}\vert}a_J(h)$ so
if $\vert \hat{I}\cap \hat{J}\vert$ is odd (knowing that $\vert
\hat{I}\vert$ is even) then $\emptyset\neq \hat{I}\cap \hat{J}\neq
\hat{I}$ and as a consequence $a_J(h)= 0$. \hfill$\square$

\begin{rem}{\rm
The most important is the case when $\hat{I}$ has two elements.
For example if $\hat{I} = \{1,2\}$ then the associated shielding
function is $f^\mu_{110\ldots 0}$. Observe that we do not need new
shielding functions to shield other singular points in the
configuration space where the action is not free. Indeed, if
$h_i=-h_j$ for some pair of indices $\{i,j\}\subset [k]$ than the
same function (\ref{eqn:shield-2}) that shields the region where
$h_i=h_j$ can be used again. For example if $i=1$ and $j=1$ then
$f^\mu_{110\ldots 0}(h) = -1$ if $h_1=-h_2$. }
\end{rem}

\subsection{Recursive computation of the parity number}
\label{sec:parity count appl}

The following result from \cite{Ram} is central tool for the
recursive computation of the parity number $P(r, X)$.

\begin{prop}{\rm (\cite[Theorem~2.2.]{Ram})}\label{prop:parity-tool}
Let $r = (r', r'') : X\rightarrow \mathbb{R}^{n-1}\times
\mathbb{R}$ be a function satisfying all the conditions listed in
Proposition~\ref{prop:standard}, including the `standard
genericity assumption'.

 \begin{enumerate}
 \item[{\rm (i)}]
Suppose that $\partial(X) = \bigcup_{i=1}^s~Y_i$ where $Y_i$ are
pairwise,  interior disjoint subcomplexes, ${\rm int}(Y_i)\cap
{\rm int}(Y_j)=\emptyset$ for $i\neq j$. Assume that ${\rm
int}(Y_i)$ are open manifolds and that $r'$ has no zeros in the
union of all boundaries $\bigcup_{i=1}^s~(Y_i\setminus{\rm
int}(Y_i))$. Then,
\begin{equation}\label{eqn:par-1}
P(r, X) = \sum_{i=1}^s P^+(r', \underline{r}''; Y_i).
\end{equation}
 \item[{\rm (ii)}] Let $Z_{r'}= r'^{-1}(0)$ be the zero-set of the function $r'$.
Suppose that $Y_i = Y_{i,1}\cup Y_{i,2}$ is an interior disjoint
union such that $r'$ has no zeros in the set
$Y_{i,\epsilon}\setminus{\rm int}(Y_{i,\epsilon})$ {\rm (}for
$\epsilon\in\{0,1\}${\rm )}. Assume that there exists a bijection
$\beta : Y_{i,1}\cap Z_{r'}\rightarrow Y_{i,2}\cap Z_{r'}$ {\rm
(}usually a restriction of a homeomorphism $\beta : Y_{i,1}
\rightarrow Y_{i,2}${\rm )} such that for some integer $a$,
$r''(\beta(x))=(-1)^ar''(x)$ for each $x\in Y_{i,1}\cap Z_{r'}$.
Then,
\begin{equation}\label{eqn:par-2}
P^+(r', \underline{r}''; Y_i) = a\cdot P(r'; Y_{i,1}).
\end{equation}
 \end{enumerate}
\end{prop}

\medskip\noindent
{\bf Proof:} Both statements are immediate consequences of
Proposition~\ref{prop:standard}. \hfill$\square$

\begin{rem}\label{rem:zastita_regiona}
{\rm In order to apply Proposition~\ref{prop:parity-tool} one
should be able to guarantee that there are no zeroes of the
function $r'$ in any of the sets $Y_i\setminus{\rm int}(Y_i)$
(respectively $Y_{i,\epsilon}\setminus{\rm int}(Y_{i,\epsilon})$).
This is usually achieved by one of the following requirements.
 \begin{enumerate}
 \item[{(1)}]The map $r'$ is generic and ${\rm dim}(Y_i\setminus{\rm int}(Y_i)) < {\rm
 dim}(Y_i)$;
 \item[{\rm (2)}] A component of $r'$ is a {\em shielding
 function}.
 \end{enumerate}
  }
\end{rem}

\begin{exam}\label{exam:BU}{\rm
As an illustration we demonstrate how the Borsuk-Ulam theorem
follows from the parity calculation described in
Proposition~\ref{prop:parity-tool}. To this end we show that if $f
: B^n\rightarrow\mathbb{R}^n$ is a (generic) map, which is
$\mathbb{Z}_2$-equivariant in the sense that $f(-x)=-f(x)$ for
each $x\in S^{n-1}$, then $P(f, B^n ) = 1$. Since the action is
free here we do not need shielding functions, i.e.\ the condition
(1) from Remark~\ref{rem:zastita_regiona} is sufficient. If $f =
(f', f'') : B^n\rightarrow \mathbb{R}^{n-1}\times \mathbb{R}$ then
by Proposition~\ref{prop:parity-tool}, $$P(f, B^n) = P^+(f',
\underline{f}''; S^{n-1}) = P^+(f', \underline{f}''; S^{n-1}_+) +
P^+(f', \underline{f}''; S_{-}^{n-1})$$ where $S^{n-1}_+$ and
$S^{n-1}_-$  are the hemispheres. Since $f''(-x) = -f''(x)$ for
$x\in S^{n-1}$ it follows from
Proposition~\ref{prop:parity-tool}(ii) that $P(f, B^n) = P(f',
S^{n-1}_+)$. Since $P(g, B^1)=1$ for a function $g :
[-1,1]\rightarrow \mathbb{R}$ satisfying the condition
$g(-1)=-g(1)$ the proof is completed by induction. }
\end{exam}

\begin{exam}\label{exam:shield}{\rm
In this example we review the proof of $\Delta(2,2) = 3$ based on
\cite[Section~5.1.]{Ram} (following the notation from this paper).
We outline some of the key steps illuminating the role of the
shielding functions.

\smallskip The inequality $\Delta(2,2)\leq 3$ is deduced from
$P(\mathbf{r}, (B^2)^2_{\leq}) = 1$ where $$(B^2)^2_{\leq} =
\{(x_1,x_2)\in \mathbb{R}^2\mid 0 \leq \|x_1 \|\leq \|x_2 \|\leq
1\}$$ and
\begin{equation}\label{eqn:ex-TM}
\mathbf{r} = (r_{11}^{\mu_1}, r_{01}^{\mu_2}, r_{10}^{\mu_2},
r_{11}^{\mu_2}) : (B^2)^2_{\leq} \rightarrow \mathbb{R}^4
\end{equation}
is the associated test-map. The calculation begins with the
observation that,
\begin{equation}\label{eqn:Ram-5-1-1}
P(\mathbf{r}, (B^2)^2_{\leq}) =   P^+(\underline{r_{11}^{\mu_1}},
r_{01}^{\mu_2}, r_{10}^{\mu_2}, r_{11}^{\mu_2}; X_{1,2}) +
P^+(\underline{r_{11}^{\mu_1}}, r_{01}^{\mu_2}, r_{10}^{\mu_2},
r_{11}^{\mu_2}; X_{2,3})
\end{equation}
where $X_{1,2}$ (respectively $X_{2,3}$) are the subsets of
$(B^2)^2_{\leq}$ satisfying the condition $\|x_1\| = \|x_2\|$
(respectively $\|x_2\| = 1$). Here we use the full power of
Remark~\ref{rem:zastita_regiona}, especially the fact that a
component $r_{11}^{\mu_2}$ of the (reduced) test map $\mathbf{r'}=
(r_{01}^{\mu_2}, r_{10}^{\mu_2}, r_{11}^{\mu_2})$ is a shielding
function. The author continues by showing that,
\begin{equation}\label{eqn:Ram-5-1-2}
P^+(\underline{r_{11}^{\mu_1}}, r_{01}^{\mu_2}, r_{10}^{\mu_2},
r_{11}^{\mu_2}; X_{2,3}) = 0
\end{equation}
which relies on the fact that for each $(x_1,x_2)\in X_{2,3}$,
$$
\mathbf{r'}(x_1,x_2) = 0 \, \Longleftrightarrow \,
\mathbf{r'}(x_2,x_1) = 0 \quad \mbox{\rm and} \quad
r_{11}^{\mu_1}(x_1,x_2) = r_{11}^{\mu_1}(x_2,x_1).
$$
Here he again uses the fact that $r_{11}^{\mu_2}$ is a shielding
function which in light of Proposition~\ref{prop:property-4.5}
rules out the possibility $\mathbf{r'}(x,x)=0$. }
\end{exam}

\begin{rem}\label{rem:justify}{\rm
It is interesting to compare two possible ways of justifying the
parity calculation in Example~\ref{exam:shield}. In our approach
we would prefer to keep the equivariance and sacrifice the
transversality of the test map on the singular set. This is always
possible as demonstrated in Section~\ref{sec:appendix}. Ramos
would prefer to preserve transversality and offer as a sacrifice
the equivariance of the test map on the singular set (or in its
vicinity). In both approaches it is the presence of the shielding
function $r_{11}^{\mu_2}$ that keeps the zeroes of $\mathbf{r'}$
away from the troublesome region. Both approaches are correct and
lead to correct calculations. }
\end{rem}

\subsection{`Counterexample' to \cite[Lemma~6.2.]{Ram}}
\label{sec:kontraprimer}

Generalizing the calculations used in his proof of $\Delta(2,2)=3$
(outlined in our Example~\ref{exam:shield}) Ramos formulated and
proved the following proposition (Lemma~6.2. in
\cite[Section~6]{Ram}). Our formulation is identical to the
original with the addition of the condition (tacitly assumed
throughout the whole of section Section~6 in \cite{Ram}) that for
each of the symmetries involved a component of the reduced test
map $\mathbf{r'}$ is a {\em shielding function} for the associated
singular set. (The meaning of the phrase {\em `symmetric for zeros
on the boundary'} is explained in \cite[page 162]{Ram}, prior to
Lemma~6.2.)

\begin{prop}{\rm (\cite[Lemma~6.2.]{Ram})}\label{prop:Ram-Lemma6.2}
Suppose that $\mathbf{r} = (\mathbf{r'}, \underline{r''}):
(B^n)^k_{\leq} \rightarrow \mathbb{R}^{nk}$ is a NPL {\em
(non-degenerate piecewise linear)} map which is {symmetric for
zeros in the boundary} and let $a\in \{0,1\}$ be the antipodality
character of $r''$ with respect to the $k$-th ball. Assume that
for each symmetry $\beta$ involved there is a component of
$\mathbf{r'}$ acting as a shielding function for the region
$S_\beta = \{x\mid \beta(x) = x\}$. Then,
\begin{equation}
P(\mathbf{r'},r''; (B^n)^k_{\leq}) = a\cdot P(\mathbf{r'};
(B^n)^{k-1}_{\leq}\times B^{n-1}).
\end{equation}
\end{prop}

\medskip\noindent
{\bf Proof:} Following \cite{Ram} the lemma is a direct
consequence of Proposition~\ref{prop:parity-tool}. \hfill
$\square$

\bigskip The key objection of \cite{BFHZ-1} to the `parity
calculation method' of Ramos was summarized and exemplified by
their `counterexample' to his \cite[Lemma~6.2]{Ram} (our
Proposition~\ref{prop:Ram-Lemma6.2}), see Example~7.7 in
\cite[p.~22]{BFHZ-1}.

\medskip
We claim that this `counterexample' is not correct in the sense
that it does not address properly Lemma~6.2. Explicitly, the map
$r = (r',r''): (B^1)^3_\leq \rightarrow \mathbb{R}^3$ they
describe does not satisfy the condition that sufficiently many
components of the reduced map $r'$ are shielding
functions!\footnote{This is evident already from the fact that the
zeroes of $r'$ are in the union of boundary sets $F_{x,y},
F_{y,z}, F_{x,z}$ which all should be shielded since they consist
solely of singular points.}

\medskip Nevertheless let us take a look at their argument more closely.
The authors of \cite{BFHZ-1} summarize their objection by saying that their example
{\em `exploits the simple fact that the permutation action on the
coordinates in $C_{m;n}$ has fixed points, a fact that Ramos does
not account for in his proof'.}

\medskip This assertion is clearly incorrect since the whole concept of a
shielding function (shield function) was invented for this
purpose. A quick computer search through \cite{BFHZ-1} reveals
that the word {\em `shield'} (as in `shield function') is
completely absent from their paper, in particular it is not
mentioned in their `counterexample' to Lemma~6.2.  It appears that
the authors of \cite{BFHZ-1} unfortunately did not read \cite{Ram}
carefully enough and apparently missed to observe the central role
played by shielding functions in this paper.

\medskip An objective reader may correctly remark that after all
there is a missing condition in \cite[Lemma~6.2.]{Ram}. We could
agree with this to some extent, however this is hardly an
`essential gap' leading to the conclusion that the `approach
employed cannot work'.

\medskip
Moreover, at the end of the paragraph (typeset in the fine print)
immediately after the proof of Lemma~6.2.\ the reader will find
the following lines:

\medskip
\noindent\hspace{1cm} \parbox[c]{14cm}{\small Thus, it is correct
to assume that the NPL (non-degenerate piecewise linear)
approximations have the required symmetry properties as long as in
the expansion in which symmetry is used, a {\em shield function
remains for each symmetry used}. }

\medskip
In other words Ramos reiterates the importance of shielding
functions and formulates precisely the `missing condition' from
his Lemma~6.2.

\section{Our response to \cite{BFHZ-1}}\label{sec:remarks}

The paper \cite{BFHZ-1} is welcome as an invitation to an
interesting and important problem in geometric combinatorics,
however it leaves much to be desired on the level of careful and
accurate presentation and interpretation of earlier proofs and
results.

\medskip
This is a pity since the criticism is always welcome, as it
provides an opportunity to improve the presentation and test one's
overall understanding of the problem.

\medskip We agree that the exposition in all papers \cite{Ram, Guide2, MVZ,
Z08} can be improved, notably in our papers \cite{Guide2, MVZ,
Z08} we tacitly used the assumption that all test maps arise from
measures (see Sections~\ref{sec:shield-principle} and
\ref{sec:gaps-corrigenda}).

\medskip

However we strongly disagree with the negative conclusions from
\cite{BFHZ-1}. We claim that the insight, basic constructions and
the results in \cite{Ram, Guide2, MVZ, Z08} are correct. The same
applies to the paper \cite{Ziv2015} (a detailed analysis of the
methods used in this paper is postponed for a subsequent paper).

\subsection{The `gaps' and `corrigenda'}\label{sec:gaps-corrigenda}

The reader may wonder how is it possible that the `essential gaps'
in so many papers passed unnoticed until the appearance of
\cite{BFHZ-1}. The answer is there are no `essential gaps' in
these papers. Here we offer a  footnote size `corrigendum'
summarizing what was said about `shielding functions' in previous
sections.

\medskip
\begin{enumerate}
 \item[{\rm (1)}] All equivariant test maps (see our
 Section~\ref{sec:set-up} and \cite[Section~4]{Ram}) arise from
 measures. As a consequence if some of the hyperplanes coincide
 some of the hyperorthants are degenerated and have measure zero.

 \item[{\rm (2)}] The test maps, restricted to
 the singular set (where some hyperplanes coincide) are therefore
 linearly homotopic (the `shielding functions homotopy principle'
 (Definition~\ref{def:shield} in Section~\ref{sec:shield-principle})).

 \item[{\rm (3)}] One uses the relative, rather than the absolute
 equivariant obstruction theory, as explicitly suggested already in our original paper
  \cite[Remark~4.3]{Guide2} (see also the introductory part of
 Section~\ref{sec:obstr-omega}).

\end{enumerate}
The reader may ask what is the main content of papers
\cite{Guide2, MVZ, Z08} (if the assumptions (1)--(3) are tacitly
treated there as part of the overall (technical) set-up). The
answer is that the real challenge in the
Gr\"{u}nbaum-Hadwiger-Ramos hyperplane mass partition problem is
always the concrete {\bf evaluation} of the topological
obstruction. Here are some highlights.

\begin{enumerate}
\item[1.] The central new idea introduced in \cite{Guide2} and
developed in \cite{MVZ, Z08} is to use measures supported by the
moment curve for the evaluation of the obstruction. Together with
the `parity count method'  of Ramos this is still the only general
method used by all papers including \cite{BFHZ-2}.

 \item[2.] The use of the moment curve (or any other {\em convex
 curve}) reduces the evaluation of the obstruction to a
 problem of enumerative combinatorics, namely to enumeration
 of combinatorial patterns related to Hamiltonian paths in hypercubes
 (Gray codes), see our Figure~\ref{fig:vamos-1} and compare it to
 Figure~2 in \cite{Guide2} or Figures~2 and 3 in \cite{MVZ}.

 \item[3.] The central fact leading to the main result in
 \cite{Z08} (Theorem~5.1) was the observation that the {\em unique balanced 4-bit Gray
 code} has an inner symmetry (Figures~2, 3, and 4 in \cite{Z08}).
\end{enumerate}

\subsection{Equivariant cobordism and shielding functions}

There is another important idea (point of view), more or less
explicit in \cite{Guide2, MVZ, Z08} (see for example Section~2.3.\
in \cite{Z08}), which also involves shielding functions and
explains to some extent how it happened that (1)--(3) (from
Section~\ref{sec:gaps-corrigenda}) were not more explicitly stated
among the assumptions in these papers.

\medskip
The idea explains how one can justify the use of open manifods,
for example (as in our Section~\ref{sec:shield-principle}) the use
of the configuration space $M_\delta = (S^5)^3_\delta :=
(S^5)^3\setminus S$ obtained by removing the singular orbits. The
illustrative case of $2$-equipartitions of a single measure in
$\mathbb{R}^2$ is presented in \cite[Section~2.3.]{Z08} as an
introduction and motivation for the more interesting (and complex)
`symmetric $4$-dimensional case' (treated later in the same
paper). The following details are extracted from
\cite[Section~2.3.]{Z08}.

\medskip
The open manifold $M_\delta = (S^2)^2_\delta$ parameterizes pairs
of distinct, oriented (affine) lines in $\mathbb{R}^2$. Note that
$M_\delta$ is an open, free $\mathbb{D}_8$-manifold where
$\mathbb{D}_8$ is the dihedral group of order $8$. Given a
(sufficiently regular) continuous measure $\mu_0$ on
$\mathbb{R}^2$, the associated solution set $\Sigma_{\mu_0}\subset
(S^2)^2_\delta$ of all equipartitions of $\mu_0$ is a compact
$1$-dimensional manifold (equipped with a free action of
$\mathbb{D}_8$). A (sufficiently generic) path $\{\mu_t\}_{t\in
[0,1]}$ (homotopy), connecting $\mu_0$ to any other (generic)
measure $\mu_1$, defines an equivariant cobordism $N_\mu$ between
the solution sets $\Sigma_{\mu_0}$ and $\Sigma_{\mu_1}$.

\begin{rem}\label{rem:obvious}{\rm
The manifold $N_\mu$ is a compact surface! This may appear obvious
and (on second thought) is obvious, however note that precisely
here we rely on the fact that the equipartitions of a given family
$\{\mu_t\}_{t\in [0,1]}$ of measures cannot escape to infinity.
Formally this is guaranteed by the existence of the corresponding
shielding function!}
\end{rem}

The proof is completed by showing that the obstruction cobordism
class $\theta = [\Sigma_{\mu_0}]$ represents the generator in the
group $\Omega_1(\mathbb{D}_8)\cong \mathbb{Z}_4$ of equivariant
bordisms. This is done by choosing the unit disc $D^2$ as the test
measure $\mu_0$ or alternatively (which is more suitable for
generalizations) its boundary $S^1$ which is an example of a {\em
convex curve} in $\mathbb{R}^2$.

\subsection{The `gaps' in the paper \cite{Ram} of Ramos}

On the basis  of the analysis presented in
Section~\ref{sec:Ramos}, see in particular
Section~\ref{sec:kontraprimer}, we conclude that there are no gaps
in the paper \cite{Ram} of Ramos. Moreover his `parity count
method' remains  a rigorous and powerful tool for proving results
about equipartitions of masses by hyperplanes.

\subsection{The role of shielding functions in \cite{BFHZ-1}}

As observed in Section~\ref{sec:kontraprimer} a computer search
through \cite{BFHZ-1} shows that the word {\em `shield'} or {\em
`shielding'} (function) is completely absent from \cite{BFHZ-1}.
It appears that the authors of \cite{BFHZ-1} completely overlooked
one of the key technical ingredients in all papers \cite{Ram,
Guide2, MVZ, Z08}. This may explain why they spent a lot of time
and energy in \cite[Section~6]{BFHZ-1} discussing
(counter)examples which have little to do with actual methods used
in \cite{Ram, Guide2, MVZ, Z08}.

\medskip
Note that these rather technical examples (see
\cite[Section~6]{BFHZ-1}) are hardly surprising to experts
interested in generalizations of the Borsuk-Ulam theorem. For
example in \cite[Theorem~6.1.]{BFHZ-1} they establish the
existence of a $\Sigma_k^\pm$-equivariant map $Z_{d,k}\rightarrow
S(U_k^{\oplus j})$ and in particular an equivariant map
$(S^4)^4_\delta \rightarrow S(U_4)$. This is more complicated but
otherwise similar in spirit to the `Borsuk-Ulam example'
introduced in Section~\ref{sec:criticism} claiming that there
exists a $\mathbb{Z}_2$-equivariant map from $S^3_\delta =
S^3\setminus\{a, -a\}$ to $\mathbb{R}^3\setminus\{0\}$.

\section{Appendix}\label{sec:appendix}

In this section we collect some basic, equivariant transversality
lemmas used throughout the paper. Our primary objective is to
provide on overview and some technical details needed for the
proof of Proposition~\ref{prop:claim} and for the applications in
Section~\ref{sec:Ramos}. The exposition is elementary and fully
accessible to a non-expert, including a combinatorially minded
reader without much previous exposure to algebraic topology.

\subsection{Equivariant transversality theorem}
\label{eqn:equi-transversality}

Let $G$ be a finite group and suppose that $M$ is a
$n$-dimensional, compact, smooth $G$-manifold. Let $S\subset M$ be
the `singular set' $S = \{x\in M \mid G_x\neq e\}$ of points with
a non-trivial stabilizer.

\medskip
Suppose that $V$ is a real, $n$-dimensional representation of $G$
and let $g : S \rightarrow V\setminus\{0\}$ be a continuous,
$G$-equivariant map.

\medskip
Suppose that $f : M\rightarrow V$ is a $G$-equivariant extension
of $g$ which are transverse to $\{0\}\in V$. The set $f^{-1}(0)$
is finite and $G$-invariant. The number of $G$-orbits $n_f = \vert
f^{-1}(0)/\vert G\vert$ clearly depends on $f$. The following
elementary lemma  claims that the parity of this number $\theta =
\theta_f =_{mod\, 2} n_f\in \mathbb{Z}_2$ is the same for all
extensions $f$ of $g$ in the same relative $G$-homotopy class.

\begin{prop}\label{prop:central}
Suppose that $f_0, f_1 : M\rightarrow V$ are two $G$-equivariant
extensions of $g$ which are transverse to $\{0\}\in V$. Then,
\begin{equation}\label{eqn:parity}
\theta(f_0) \equiv \theta(f_1) \qquad \mbox{ {\rm mod} } (2).
\end{equation}
\end{prop}

\medskip\noindent
{\bf Proof:} Let $F: M\times [0,1] \rightarrow V$ be the linear
homotopy between $f_0$ and $f_1$ defined by $F(x,t) = (1-t)f_0(x)
+ t f_1(x)$. Let $A = (S\times [0,1])\cup (M\times\{0,1\})$. The
map $F$ is by assumption already transverse to $\{0\}\in V$ on the
set $A$. Since the action of the group $G$ is free on $(M\times
I)\setminus (S\times I)$, we are allowed to apply the `Equivariant
Transversality Theorem' (Theorem~\ref{thm:ETT}) which claims that
the homotopy $F$ admits a small $G$-equivariant perturbation,
\begin{equation}
H : M\times [0,1] \rightarrow V
\end{equation}
which is transverse to $0\in V$ everywhere. Moreover we can assume
that $H$ and $F$ agree on the set $A$. It follows that $H^{-1}(0)$
is a compact, $1$-dimensional manifold ($G$-bordism), connecting
zero sets $Z(f_0)$ and $Z(f_1)$. In turn there is a bordism
between the associated sets $Z(f_0)/G, Z(f_1)/G$ of orbits and the
equality (\ref{eqn:parity}) is an immediate consequence. \hfill
$\square$

\begin{theo}{\rm(Equivariant Transversality Theorem)}\label{thm:ETT}
Let $G$ be a finite group and suppose that $N$ is a
$n$-dimensional, compact, smooth $G$-manifold. Let $S\subset N$ be
the `singular set' $S = \{x\in N \mid G_x\neq e\}$ of points with
a non-trivial stabilizer. Suppose that $V$ is a real,
$k$-dimensional representation of $G$. Suppose that,
\[
F : N \rightarrow V
\]
is a $G$-map such that $0\notin F(S)$. Suppose that $A\subset N$
is a closed $G$-subset where $F$ is already transverse to $\{0\}$.
Then there exists a $G$-equivariant map,
\[
H : N \rightarrow V
\]
transverse to $\{0\}$ which agrees with $F$ on $A\cup S$.
Moreover, there exists a relative $G$-homotopy (small perturbation
of $F$),
\[
  G : N\times [0,1] \rightarrow V \qquad \mbox{ {\rm (rel} } A\cup S)
\]
connecting $F$ and $H$.
\end{theo}

\medskip\noindent
{\bf Comments on the proof of Theorem~\ref{thm:ETT}:}
Theorem~\ref{thm:ETT} is quite directly a consequence of  the
standard (non-equivariant) transversality theorem (or its proof),
as exposed in \cite{GG, Ko} and other textbooks.

\medskip
One of the guiding principles used in the standard proofs (see for
example \cite{GG}) is to make the map transverse locally (one
small open set at a time) using the fact that a small perturbation
will not affect the transversality condition achieved earlier on
in the construction.

This can be done equivariantly, in the region where the action is
{\em free}. Indeed, if $V\subset N$ is a (small) open set such
that $V\cap g(V)$ for each $g\neq e$ then $F$ can be made
transverse to $\{0\}$ on $V$ and extended equivariantly to
$\cup_{g\in G}~g(V)$, etc.

\medskip
Another possibility is to use the fact that equivariant maps are
sections of a bundle. Let $U$ be the interior of a $G$-invariant
regular neighborhood of $S$. Let $N' = N\setminus U$ and $\partial
N' = \partial \overline{U}$. We can assume that $0\notin
F(\partial N')$. Then the action of $G$ on $N'$ is free and there
is a one-to-one correspondence between $G$-equivariant maps $f :
N'\rightarrow V$ and sections of the bundle,
\begin{equation}\label{eqn:bundle}
        V \longrightarrow N'\times_G V \rightarrow N'/G .
\end{equation}
Moreover $f$ is transverse to $\{0\}$ if and only if $s$ is
transverse (in the usual, non-equivariant sense) to the zero
section of the bundle. \hfill $\square$

\end{document}